\documentstyle [12pt]{article}
\begin{document}
\title{
The space of free coherent states
\\
is isomorphic to the space
\\
of distributions on $p$-adic numbers
}
\author{S.V.Kozyrev}
\maketitle

\footnotetext{
e-mail: kozyrev@genesis.mian.su,\quad Address:  Inst.Chem.Phys., OSV, 117334,
Kosygina 4, Moscow}

\begin{abstract}
Free coherent states for a system with $p$ degrees of freedom are defined.
An isomorphism  of the space of free coherent states to the space of
distributions on $p$-adic disk  is constructed.
\end{abstract}

\section{Construction of free (or Boltzmannian) coherent states}

Free (or Boltzmannian)
Fock space has been considered in some recent works
on quantum chromodynamics
\cite{MasterFld}, \cite{GopGross}, \cite{DougLi}
and noncommutative probability
\cite{AccLu}, \cite{Maassen}, \cite{Speicher}.

The subject of this work is free coherent states.
We will introduce free coherent states and investigate
the  space of coherent states corresponding to a fixed eigenvalue
of the operator of annihilation.

In this paper the connection between the theory of free Fock space and
$p$-adic analysis will be established.
The main result of the present paper is a construction of
an isomorphism of the space of coherent states to
the space of distributions on the ring of integer $p$-adic numbers.

We will consider the system with $p$ degrees of freedom.
The system with one degree of freedom was investigated in \cite{KozAAV}.

The free commutation relations are particular case of
$q$-deformed relations
$$A_i A_j^{\dag}-qA_j^{\dag}A_i=\delta_{ij}  $$
with $q=0$.
A correspondence of  $q$-deformed commutation relations and
non-archimedean (ultrametric) geometry was discussed in
\cite{ArVqGPnAG}.
Non-archimedean mathematical physics
was studied in \cite{VVZ}.

Free coherent states lie in
the free Fock space. Free  (or Boltzmannian) Fock space
$F$ over a Hilbert space  $H$ is the completion of the tensor algebra
$$F=\oplus_{n=0}^{\infty}H^{\otimes n}.$$
Creation and annihilation operators are defined in the following way:
$$
A^{\dag}(f) f_{1}\otimes \dots \otimes f_{n}=f\otimes
f_{1}\otimes \dots \otimes f_{n}
$$
$$
A(f) f_{1}\otimes \dots \otimes
f_{n}=\langle f,f_{1} \rangle f_{2}\otimes \dots  \otimes f_{n}
$$
where  $\langle f,g \rangle $ is the scalar product in the Hilbert space $H$.
Scalar product  in the free Fock space is defined by the standard
construction of the direct sum of tensor products of Euclidean spaces.

We consider the case when  $H$ is $p$-dimensional Euclidean space.
In this case we have $p$ creation operators
$A^{\dag}_{i}$, $i=0,\dots ,p-1$
and $p$ annihilation operators
$A_{i}$, $i=0,\dots ,p-1$ with commutation relations
\begin{equation}\label{aac}
A_{i}A_{j}^{\dag}=\delta_{ij}.
\end{equation}
The vacuum vector $\Omega$  in the free Fock space satisfies
\begin{equation}\label{vacuum}
A_{i}\Omega=0.
\end{equation}

We introduce
the linear space   of free coherent states
as the space of formal eigenvectors of the
annihilation operator $A=\sum_{i=0}^{p-1}A_{i}$ for the eigenvalue $\lambda$,
$$
A \Psi= \lambda \Psi.
$$
It is easy to see that the formal solution of this equation is
\begin{equation}\label{psi}
\Psi=\sum_{I} \lambda^{|I|} \Psi_I A^{\dag}_{I}\Omega.
\end{equation}
Here  the multiindex $I=i_0 \dots i_{k-1}$, $i_j \in\{0, \dots ,p-1\}$ and
$A^{\dag}_{I}=A^{\dag}_{i_{k-1}} \dots A^{\dag}_{i_0}$,
$\Psi_I$ are complex numbers that satisfy
$$
\Psi_I=\sum_{i=0}^{p-1}\Psi_{Ii}.
$$
Summation in the formula  (\ref{psi}) runs on all sequences $I$ with the
finite length. The length of a sequence $I$ is denoted by $|I|$.
This formal series define a functional with
a dense domain in the free Fock space.

An arbitrary free coherent state $\Psi$ can be constructed in the following
way. $\Psi$ is defined by the function  $\Psi_{I}$.
Let us construct this function.
Let us define $\Psi_{\emptyset}$  as an arbitrary complex number.
Then we use the inductive procedure: if we know $\Psi_{I}$
we define $\Psi_{Ii}$ as  arbitrary complex numbers
that satisfy
the formula $\Psi_{I}=\sum_{i=0}^{p-1}\Psi_{Ii}$.
It is easy to see that the function  $\Psi_{I}$
satisfies the following formula
\begin{equation}\label{cascade}
\Psi_{I}=\sum_{|J|=j}\Psi_{IJ}
\end{equation}
for arbitrary $j$.

We will denote the linear space of free coherent states as $X'$
because this space is isomorphic to the space $D'(Z_p)$ of  distributions
on the $p$-adic disk $Z_p$.

\section{A brief review of $p$-adic analysis}

Let us make a brief review of  $p$-adic analysis.
The field $Q_p$ of $p$-adic numbers is the completion of the field of rational
numbers  $Q$ with respect to the $p$-adic norm on $Q$.
This norm is defined in the following way. An arbitrary rational number
$x$ can be written in the form $x=p^{\gamma}\frac{m}{n}$ with $m$ and $n$
that are not divisible by $p$. The $p$-adic norm of the rational number
$x=p^{\gamma}\frac{m}{n}$ is equal to $||x||_p=p^{-\gamma}$.

The most interesting property of the field of   $p$-adic numbers is
ultrametricity. This means that $Q_p$ obeys the strong triangle inequality
$$
||x-y||_p \le \max (||x-z||_p,||z-y||_p).
$$
We will consider disks in   $Q_p$ of the form
$\{x\in Q_p: ||x-x_0||_p\le p^{-k}\}$.
For example, the ring $Z_p$ of integer $p$-adic numbers
is the disk
$\{x\in Q_p: ||x||_p\le 1\}$.
The main properties of disks in
arbitrary ultrametric space are the following:

{\bf 1.}\qquad
Every point of a disk is the center of this disk.

{\bf 2.}\qquad
Arbitrary two disks either do not intersect or one of these disks
contains  another.

The main result of the present paper is the isomorphism of
the space $X'$ of free coherent states and
the space $D'(Z_p)$
of distributions on the $p$-adic disk.
The space $D'(Z_p)$
of distributions on the $p$-adic disk
is the space of linear functionals on the linear space of locally constant
complex valued functions.
The space of locally constant functions
is the linear span of indicators of $p$-adic disks with the radius
that less  or equal to 1.
Here indicator of a set $S$  is the function that equals to 1 on this set
and equals to 0 outside the set $S$.

For further reading on the subject of  $p$-adic analysis see
\cite{Vlad}, \cite{VVZ}.

\section{Construction of an isomorphism of  the space
of coherent states to the space of distributions on the $p$-adic disk}

In the present section the isomorphism $\phi$ of the space $X'$
of free coherent states to the space of distributions on the $p$-adic disk
$Z_p$ of integer $p$-adic numbers
will be constructed. In the present section we take 
the eigenvalue $\lambda$, see (\ref{psi}),  
from the interval $(0,\sqrt{p})$.

Distributions on the $p$-adic disk are linear functionals on the space
of locally constant functions \cite{Vlad}, \cite{VVZ}.
Therefore we have to construct the coherent state that corresponds to a locally
constant function.
It is sufficient to find coherent states that correspond to locally constant
functions of the type
$$
\theta_k(x-x_0)=\theta(p^{k}||x-x_0||_p);\quad
\theta(t)=0, t>1;\quad \theta(t)=1, t\le 1.
$$
Here $x$, $x_0\in Z_p$ lie in
the ring of integer $p$-adic numbers and the function $\theta_k(x-x_0)$ equals
to 1 on the disk $D(x_0,p^{-k})$ of radius $p^{-k}$ with the center in $x_0$
and equals to 0 outside this disk.

Let us consider the multiindex $I=i_0 \dots i_{k-1}$, $i_j=0, \dots ,p-1$.
Let us introduce the free coherent state $X_I$ of the form

\begin{equation}\label{indicator}
X_I=
\sum_{k=0}^{\infty} \lambda^k
\left(\frac{1}{p} \sum_{i=0}^{p-1}A_i^{\dag}\right)^k
\lambda^{|I|} A^{\dag}_I \Omega+
\sum_{l=1}^{\infty} \lambda^{-l}
\left(\sum_{i=0}^{p-1}A_i\right)^l
\lambda^{|I|} A^{\dag}_I \Omega.
\end{equation}
The sum on $l$ in fact contains $|I|$ terms.
Under the isomorphism $\phi$ mentioned above the coherent state
$X_I$ will correspond to the locally constant function
$\frac{\theta_{|I|}(x-I)}{||\theta_{|I|}(x-I)||^{{2}}}$,
where
$||\theta_{|I|}(x-I)||$ is  $L_2$-norm.
Here we identify the sequence
$I=i_0 \dots i_{k-1}$  and the $p$-adic number  $I=\sum_{j=0}^{k-1}i_j p^j$.

For $\lambda\in (0,\sqrt{p})$ the coherent state $X_I$ lies in the
Hilbert space (the correspondent functional is bounded).

We will prove that
an arbitrary locally constant function  corresponds to the
linear combination of coherent
states $X_I$.
The linear span of
the coherent states $X_I$
we will denote by $X$.

Every vector in  $X$ is a function of $\lambda$. We will investigate
the  properties of the space $X$ with the scalar product
\begin{equation}\label{product_trun}
\langle X_I,X_J \rangle=
\lim_{\lambda\to\sqrt{p}-0}\left(1-\frac{\lambda^2}{p}\right)
(X_I,X_J).
\end{equation}
Here $(X_I,X_J)$ is the scalar product in the free Fock space.

We will prove the following lemma.

{\bf Lemma 1.}\qquad
{\sl
The  limit of
the scalar product of coherent states $X_I$, $X_J$
equals to the integral on $p$-adic disk with respect to the Haar measure
$$
\lim_{\lambda\to\sqrt{p}-0}\left(1-\frac{\lambda^2}{p}\right)
(X_I,X_J)=
$$
$$
=\frac{1}{\mu(D(I,p^{-|I|}))\mu(D(J,p^{-|J|}))}
\int_{Z_p}\theta_{|I|}(x-I)\theta_{|J|}(x-J)dx=
$$
$$
=\left(
\frac{\theta_{|I|}(x-I)}{||\theta_{|I|}(x-I)||^{{2}}},
\frac{\theta_{|J|}(x-J)}{||\theta_{|J|}(x-J)||^{{2}}}
\right)_{L_2}.
$$
Here $\mu(D)$ is the Haar measure of the disk $D$.

}

{\bf Proof}

Let us calculate the scalar product of coherent states
$X_I$ and $X_J$.
Let $|I|\le |J|$.
If the  the sequence $I$  coincides with
the first $|I|$ indices of the sequence $J$ then the scalar product
has the following form
$$
(X_I,X_J)=\sum_{i=0}^{|I|-1}\lambda^{2i}+
\sum_{i=|I|}^{|J|-1}\lambda^{2i}\left(\frac{1}{p}\right)^{i-|I|}+
\sum_{i=|J|}^{\infty}\lambda^{2i}\left(\frac{1}{p}\right)^{i-|I|}.
$$
If the  sequence $I$  does not coincide with
the first $|I|$ indices of the sequence $J$ then the series for scalar product
$(X_I,X_J)$ contains only the finite number of terms.
Therefore the limit
$\lim_{\lambda\to\sqrt{p}-0}\left(1-\frac{\lambda^2}{p}\right)
(X_I,X_J)$
equals to $\min(p^{|I|},p^{|J|})$ if one of the disks $D(I,p^{-|I|})$
and $D(J,p^{-|J|})$  contains another and equals to 0 if these disks
do not intersect.

Therefore the  limit of the scalar product of vectors $X_I$,
$X_J$ in the free Fock space equals to the integral on $p$-adic disk
with respect to the Haar measure. Coherent state $X_I$ corresponds to
the locally constant function
$\frac{1}{\mu(D(I,p^{-|I|}))}\theta_{|I|}(x-I)=
\frac{\theta_{|I|}(x-I)}{||\theta_{|I|}(x-I)||^{{2}}}
$.

Let us investigate functionals on the space $X$ of locally constant functions.
Let us consider an infinite sequence $I=i_0 \dots i_{k}...$,
$i_j=0, \dots ,p-1$
that corresponds to the $p$-adic number $I=\sum_{k=0}^{\infty}i_k p^k$.
Let us denote $I_k=i_0 \dots i_{k-1}$.
Let us introduce the free coherent state $X_I$ of the form
$$
X_I= \sum_{k=0}^{\infty}
\lambda^{|I_k|} A^{\dag}_{I_k} \Omega.
$$
Let us consider the coherent state $X_J$, $J=j_0 \dots j_{|J|-1}$.

{\bf Lemma 2.}\qquad
{\sl
The limit of  the action of the functional $X_I$
on the vector  $X_J$
has the following form
$$
\lim_{\lambda\to\sqrt{p}-0}\left(1-\frac{\lambda^2}{p}\right)
(X_I,X_J)=
\frac{1}{\mu(D(J,p^{-|J|}))}\int_{Z_p}\delta(x-I)\theta_{|J|}(x-J)dx.
$$
}

Therefore the coherent state $X_I$ corresponds to
the $\delta$-function   $\delta(x-I)$.

The next lemma allows us to identify the space of free coherent states
and the space of distributions on a $p$-adic disk.

{\bf Lemma 3.}\qquad
{\sl
Vectors $X_I\in X$ lie in the domain of the functional $\Psi$
for an arbitrary free coherent state $\Psi$.
}

{\bf Proof}

The  coherent state
$X_I$ has the following form
$$
X_I=
\sum_{k=0}^{\infty} \lambda^k
\left(\frac{1}{p} \sum_{i=0}^{p-1}A_i^{\dag}\right)^k
\lambda^{|I|} A^{\dag}_I \Omega+
\sum_{l=1}^{\infty} \lambda^{-l}
\left(\sum_{i=0}^{p-1}A_i\right)^l
\lambda^{|I|} A^{\dag}_I \Omega.
$$
The action of the functional $\Psi$
is defined by the following formal series
\begin{equation}\label{action}
(\Psi,X_I)=\sum_{k=0}^{\infty}\lambda^{2k} (\Psi^{k},X_I^{k}).
\end{equation}
Here $\Psi^{k}$ and $X_I^{k}$ are coefficients of $\lambda^k$
in series for $\Psi$ and $X_I$.
$\Psi^{k}$ is defined by the formula
$$
\Psi^{k}=\sum_{|J|=k}  \Psi_J A^{\dag}_{J}\Omega;
$$
and $X_I^{k}$ for $k>|I|$ have a form
$$
X_I^k=
\left(\frac{1}{p} \sum_{i=0}^{p-1}A_i^{\dag}\right)^{k-|I|}
A^{\dag}_I \Omega.
$$
Let us calculate
$$
(\Psi^{|I|+k},X_I^{|I|+k})=
(X_I^{|I|+k},\Psi^{|I|+k})^{*}=
(X_I^{|I|+k-1},\frac{1}{p}\sum_{i=0}^{p-1}A_i\Psi^{|I|+k})^{*}=
$$
$$
=(X_I^{|I|+k-1},\frac{1}{p}\sum_{i=0}^{p-1}A_i
\sum_{|J|=|I|+k}  \Psi_J A^{\dag}_{J}\Omega)^{*}
$$
for $k>0$. We have
$$
\frac{1}{p}\sum_{i=0}^{p-1}A_i
\sum_{|J|=|I|+k}  \Psi_J A^{\dag}_{J}\Omega=
\sum_{|J|=|I|+k-1}\frac{1}{p} \sum_{i=0}^{p-1} \Psi_{Ji} A^{\dag}_{I}\Omega=
\frac{1}{p}\Psi^{|I|+k-1}.
$$
Therefore
$$
(\Psi^{|I|+k},X_I^{|I|+k})=\frac{1}{p}(\Psi^{|I|+k-1},X_I^{|I|+k-1}).
$$
We have $\frac{\lambda^2}{p}<1$
and the series (\ref{action}) converges.

Therefore an arbitrary coherent state $\Psi$ corresponds to a distribution
on the $p$-adic disk.
We get the following theorem.

{\bf Theorem.}\qquad
{\sl
The map $\phi$
$$
\phi:\quad X\to D(Z_p);
$$
$$
X_I\mapsto \frac{1}{\mu(D(I,p^{-|I|}))}\theta_{|I|}(x-I);
$$
extends to the isomorphism of the space $X$  of coherent states
of a type (\ref{indicator})
onto the space $D(Z_p)$ of locally constant functions on the ring of integer
$p$-adic numbers with the scalar product defined
with respect to the Haar measure.

The scalar product of locally constant functions  in $L_2$
with respect to the Haar measure
equals to the limit of the scalar product
in the free Fock space
of corresponding
coherent states
$$
\lim_{\lambda\to\sqrt{p}-0}\left(1-\frac{\lambda^2}{p}\right)
(X_I,X_J)=
$$
$$
=\frac{1}{\mu(D(I,2^{-|I|}))\mu(D(J,2^{-|J|}))}
\int_{Z_p}\theta_{|I|}(x-I)\theta_{|J|}(x-J)dx.
$$

The map $\phi$ defines an isomorphism
of the space of free coherent states onto the space $D'(Z_p)$
of distributions on the $p$-adic disk:
$$
\phi:\quad X'\to D'(Z_p).
$$
The coherent state $\Psi$
corresponds to a distribution on the ring of integer 2-adic numbers
with the action on locally constant functions defined by the formula
$$
(\phi(\Psi),\phi(X_I))=
\lim_{\lambda\to\sqrt{p}-0}\left(1-\frac{\lambda^2}{p}\right)
(\Psi,X_I).
$$

}

{\bf Proof}

We have to prove that an arbitrary distribution on $Z_p$ has a form
$\phi(\Psi)$.
A    distribution $f$ on $Z_p$ is defined by the action of $f$
on localy constant functions
$\theta_{|I|}(x-I)$ that correspond to
coherent states ${\mu(D(I,p^{-|I|}))}X_I$. We have
$$
\lim_{\lambda\to\sqrt{p}-0}\left(1-\frac{\lambda^2}{p}\right)
(\Psi,{\mu(D(I,p^{-|I|}))}X_I) = \Psi_I.
$$
Here $\Psi_I$ can be  an arbitrary function that satisfies the formula
(\ref{cascade}).
Therefore $\phi(\Psi)$ can be an arbitrary distribution on the $p$-adic disk.

The formula (\ref{cascade}) corresponds to the
linearity of the action of the distribution $\phi(\Psi)$ on the sum
of indicators
of $p$-adic disks.

The theorem  is proved.

\vspace{3mm}
{\bf Acknowledgments}

Author is grateful to I.V.Volovich for discussions.

\end{document}